\documentclass[11pt]{article}

\usepackage{xypic}
\usepackage{amsgen}
\usepackage{amsmath}
\usepackage{amstext}
\usepackage{amsbsy}
\usepackage{amsopn}
\usepackage{amsfonts}
\usepackage{amssymb}
\usepackage{graphicx}
\usepackage{pstricks}
\xyoption{all}

\def\Box{\square}
\def\edge{\relbar\joinrel\relbar}
\def\mapright#1{\smash{\mathop{\longrightarrow}\limits^{#1}}}
\def\tra#1{\smash{\mathop{\mid\kern
-1pt\joinrel\relbar\joinrel\relbar}\limits^{*}_{#1}}}
\def\longtra#1{\smash{\mathop{\mid\kern
-1pt\joinrel\relbar\joinrel\relbar\joinrel\relbar}\limits^{*}_{#1}}}
\def\vlongtra#1{\smash{\mathop{\mid\kern
-1pt\joinrel\relbar\joinrel\relbar\joinrel\relbar\joinrel\relbar}\limits^{*}_{#1}}}
\def\vvlongtra#1{\smash{\mathop{\mid\kern
-1pt\joinrel\relbar\joinrel\relbar\joinrel\relbar\joinrel\relbar\joinrel\relbar}\limits^{*}_{#1}}}
\def\vvvlongtra#1{\smash{\mathop{\mid\kern
-1pt\joinrel\relbar\joinrel\relbar\joinrel\relbar\joinrel\relbar\joinrel\relbar\joinrel\relbar}\limits^{*}_{#1}}}
\def\etra#1{\smash{\mathop{\mid\kern
-1pt\joinrel\relbar\joinrel\relbar}\limits_{#1}}}

\def\iff{\Leftrightarrow}
\def\Rw{\Rightarrow}
\def\oo{\overline}

\def\B{{\cal{B}}}

\newcommand{\N}{{\rm I}\kern-2pt {\rm N}}

\def\cay{\mbox{Cay}}

\def\endo{\mbox{End}\,}

\def\fix{\mbox{Fix}\,}
\def\per{\mbox{Per}\,}

\def\pref{\mbox{Pref}}

\def\max{\mbox{max}}

\def\sup{\mbox{sup}}

\def\p{\varphi}

\def\inv{^{-1}}

\def\bi{\begin{itemize}}
\def\ei{\end{itemize}}
\def\beq{\begin{equation}}
\def\eeq{\end{equation}}

\xyoption{all}

\newtheorem{T}{Theorem}[section]
\newcommand{\bt}{\begin{T}}
\newcommand{\et}{\end{T}}
\newcommand{\ftd}{$\square$\end{T}}

\newtheorem{Proposition}[T]{Proposition}
\newcommand{\bp}{\begin{Proposition}}
\newcommand{\ep}{\end{Proposition}}
\newcommand{\fpd}{$\square$\end{Proposition}}

\newtheorem{Lemma}[T]{Lemma}
\newcommand{\bl}{\begin{Lemma}}
\newcommand{\el}{\end{Lemma}}
\newcommand{\fld}{$\square$\end{Lemma}}

\newtheorem{Corol}[T]{Corollary}
\newcommand{\bc}{\begin{Corol}}
\newcommand{\ec}{\end{Corol}}
\newcommand{\fcd}{$\square$\end{Corol}}

\newtheorem{Result}[T]{Result}
\newcommand{\br}{\begin{Result}}
\newcommand{\er}{\end{Result}}
\newcommand{\frd}{$\square$\end{Result}}

\newtheorem{Example}[T]{Example}
\newcommand{\be}{\begin{Example}}
\newcommand{\ee}{\end{Example}}

\newtheorem{Problem}[T]{Problem}
\newcommand{\bq}{\begin{Problem}}
\newcommand{\eq}{\end{Problem}}

\newcommand{\proof}
   {\par\medbreak\noindent{\bf Proof}.\enspace}

\newcommand{\qed}{
$\Box$
\par\bigbreak}

\normalbaselines
\def\abstract#1{\par\bigskip
\begingroup\small
\baselineskip=12truept
\begin{center}ABSTRACT\end{center}
\par\medskip\par\noindent
\null\hfill\hbox{\vbox{\hsize=5truein\noindent#1}}
\hfill\null\par\endgroup\par}

\title{Fixed points of endomorphisms of trace monoids}
\author{{\bf Emanuele Rodaro, Pedro V. Silva}\\ $ $\\ {\em Centro de
Matem\'{a}tica, Faculdade de Ci\^{e}ncias, Universidade do
Porto,}\\ {\em R. Campo Alegre 687, 4169-007 Porto, Portugal}\\
{\em e-mail:} emanuele.rodaro@fc.up.pt, pvsilva@fc.up.pt}
\date{\today}

\begin{document}
\maketitle

\begin{center}\small
2010 Mathematics Subject Classification: 20M35, 68Q85, 54E50

\medskip

Keywords: trace monoid, real traces, endomorphism, fixed point,
periodic point
\end{center}

\abstract{It is proved that the fixed point submonoid and the periodic
point submonoid of a trace monoid endomorphism are always finitely
generated. Considering the Foata normal form metric on trace monoids and
uniformly continuous endomorphisms, a finiteness theorem is proved
for the infinite fixed points of the continuous extension to real traces.}

\section{Introduction}

In \cite{CS}, Cassaigne and the second author studied finiteness
conditions for the infinite fixed points of (uniformly continuous)
endomorphisms of monoids defined by special confluent rewriting
systems, extending results known for free monoids \cite{PP}. This line
of reasearch was pursued by the second author in subsequent papers
\cite{Sil2,Sil3}.

Similar problems were considered by the second author for virtually
free groups in \cite{Sil} (see also \cite{Syk}).  

Recently, in \cite{RSS}, Sykiotis and the authors showed that, given a graph
group $G$, the subgroup of
fixed points is finitely generated for every endomorphism of $G$
if and only if $G$ is a free product of free abelian groups.

Motivated by these papers, we perform here a similar study for trace
monoids, both at finite and infinite level. We remark that trace
monoids are one of the most important 
models for parallel computing in theoretical computer science
\cite{Die,DR}. In our study, the infinite fixed points are taken
among real traces, the completion of a trace monoid for the Foata
normal form metric, introduzed by Bonizzoni, Mauri and Pighizzini
\cite{BMP}.

\section{Preliminaries}

Given a monoid $M$, we denote by $\endo M$ the endomorphism monoid of
$M$. 
Given $\p \in \endo M$, we say that $x \in M$ is a {\em fixed point}
of $\p$ if $x\p = x$. If $x\p^n = x$ for some $n \geq 1$, we say that
$x$ is a {\em periodic point} of $\p$. Let
$\fix\p$ (respectively $\per\p$) denote the set of all fixed points
(respectively periodic points) of $\p$. Clearly,
$$\per\p = \displaystyle\bigcup_{n \geq 1} \fix\p^n.$$

Given $X \subseteq M$, the {\em star} operator $X \mapsto X^*$ defines
the submonoid of $M$ generated by $X$. We say that $X \subseteq M$ is
{\em rational} if $X$ 
can be obtained from finite subsets of
$M$ by applying finitely many times the operators union, product and
star. For more details on rational subsets, namely connections with
finite automata, see \cite{Ber}.

Given an alphabet $A$, the {\em free monoid} on $A$ is denoted as
usual by $A^*$. A (finite) {\em independence alphabet} is an ordered
pair of the form $(A,I)$, 
where $A$ is a (finite) set and $I$ is a symmetric anti-reflexive relation on
$A$. We can view $(A,I)$ as an undirected graph without loops or
multiple edges, denoted by $\Gamma(A,I)$, by taking $A$ as the vertex
set and $I$ as the edge 
set. Conversely, every such graph determines an independence alphabet.

Let $\rho(A,I)$ denote the congruence on $A^*$ generated by the the relation
\beq
\label{sexta}
\{ (ab,ba) \mid (a,b) \in I \}.
\eeq
The {\em trace monoid} $M(A,I)$ is the quotient $A^*/\rho(A,I)$, i.e the
monoid defined by the monoid presentation 
$$\langle A \mid ab = ba \; ((a,b) \in I) \rangle.$$
Such monoids are also known as {\em
  free partially commutative monoids}. Fore details on trace monoids,
the reader is referred to \cite{Die,DR}.

The elements of $M(A,I)$ are often described through the {\em Foata
  normal form} (FNF), defined as follows. We say that $B \subseteq A$ is
a {\em clique} if $B \neq \emptyset$ and the full subgraph of
$\Gamma(A,I)$ induced by $B$ is 
complete. Let then $w_B \in M(A,I)$ be the product of all the letters
in $B$. Note that, since the letters of $B$ commmute with each other
in $M(A,I)$, we do not need to specify any multiplication sequence. A
FNF is either the empty word or
a product $w_{B_1}\ldots w_{B_k}$ such
  that:
\bi
\item
$B_i$ is a clique for $i = 1,\ldots,k$;
\item
for all $i = 2,\ldots,k$ and $a \in B_i$, there exists some $b \in
B_{i-1}$ such that $(a,b) \notin I$.
\ei
Then each element of $M(A,I)$ admits a unique representation in FNF \cite{CF}. 

Given $u \in M(A,I)$ and $a \in A$, we denote by
$|u|_a$ the number of occurrences of $a$ in $u$.
Given $B \subseteq A$,
 let $\pi_B \in \endo M$ keep the letters
of $B$ and remove the others.
A well-known result (see \cite{DR}) states that, for all $u,v \in
A^*$,
\beq
\label{quinta}
u = v \mbox{ in } M(A,I) \hspace{.3cm} \iff \hspace{.3cm} \forall
(a,b) \in (A\times A) \setminus I \hspace{.3cm} u\pi_{a,b} = v\pi_{a,b}.
\eeq
This leads to an embedding of $M(A,I)$ into a direct product of free
monoids. It follows that $M(A,I)$ is a
cancellative monoid.

\section{The submonoids of fixed and periodic points}

\bt
\label{fixtrace}
Let $(A,I)$ be a finite independence alphabet. Then 
{\rm Fix}$\,\p$ is finitely generated for every $\p \in$
  {\rm End}$\, M(A,I)$.
\et

\proof
Write $M = M(A,I)$.
We use induction on $|A|$. The case $|A| = 0$ being trivial, assume
that $|A| > 0$ and the result holds for smaller alphabets. We consider
three cases.

\medskip
\noindent
\underline{Case I}: $\p|_A$ is a permutation.

\medskip
\noindent
We start showing that 
\beq
\label{fixtrace2}
(a,b) \in I \iff (a\p,b\p) \in I
\eeq
holds for all $a,b \in A$. Indeed, the direct implication follows from
$\p$ being an endomorphism since $ab = ba$ yields $(a\p)(b\p) =
(b\p)(a\p)$ in $M$. The converse implication follows from
the equality
$$\begin{array}{l}
|\{ (a,b) \in A \times A \mid ab = ba\mbox{ in }M \}|\\
\hspace{2cm} =
|\{ (a,b) \in A \times A \mid (a\p)(b\p) =
(b\p)(a\p) \mbox{ in }M \}|
\end{array}$$
combined with the direct implication. Thus (\ref{fixtrace2}) holds.

Let $\B$ be the set of all cliques $B$ such that $B\p = B$. We claim
that
\beq
\label{fixtrace1}
\fix\p = \{ w_B \mid B \in \B \}^*.
\eeq
It is immediate that $w_B \in \fix\p$ for every $B \in
\B$. Hence $\{ w_B \mid B \in \B \}^* \subseteq \fix\p$. 
Conversely, let $u \in \fix\p$. Assuming that $u \neq 1$, let $u = 
w_{B_1}\ldots w_{B_k}$ be its FNF.  Then $u\p =
(w_{B_1}\p)\ldots (w_{B_k}\p)$. By (\ref{fixtrace2}),
$(w_{B_1}\p)\ldots (w_{B_k}\p)$ must be a FNF as
well. Since $u\p = u$, and the representation in FNF
is unique, we get $w_{B_i} = w_{B_i}\p$ for $i =
1,\ldots,k$. Hence $B_1,\ldots,B_k \in \B$ and so $u \in \{ w_B \mid B
\in \B \}^*$. Thus (\ref{fixtrace1}) holds and so $\fix\p$ is finitely
generated in this case.

\medskip
\noindent
\underline{Case II}: $1 \notin A\p$.

\medskip
\noindent
Write
$$A_0 = A \cap \per\p$$
and let $\Lambda$ be the directed graph with vertices in $A_0$ and edges
$a \mapright{} b$ whenever $b = a\p$. It is easy to see that $A_0$
consists of all the vertices occurring in some cycle of $\Lambda$. Let
$p = |A_0|!$. Since $p$ is a multiple of the length of any cycle in
$\Lambda$, it follows that
\beq
\label{fixtrace3}
A_0 = A \cap \fix\p^p.
\eeq
We claim next that
\beq
\label{fixtrace4}
A_0 = A \cap A\p^p.
\eeq
Indeed, $A_0 \subseteq A \cap A\p^p$ follows from
(\ref{fixtrace3}). Conversely, let $a \in A \cap A\p^p$. Since $1
\notin A\p$, there exists some path $b \mapright{} a$ of length $p$
in $\Lambda$ and so $a$ lies in some cycle of $\Lambda$. Thus $a \in
A_0$ and so (\ref{fixtrace4}) holds.

Next we show that
\beq
\label{fixtrace5}
A_0\p = A_0.
\eeq
Let $a \in A_0$. By (\ref{fixtrace4}), we have $a\p = (a\p)\p^p$ and
so $a\p \in A_0$. Hence $A_0\p \subseteq A_0$. But now 
$a\p^{p-1} \in A_0$ and  $a = (a\p^{p-1})\p$ yield $A_0 \subseteq
A_0\p$ and so (\ref{fixtrace5}) holds.

Let $M'$ be the (trace) submonoid of $M$ generated by $A_0$. In view of
(\ref{fixtrace5}), $\p$ restricts to some endomorphism $\p'$ of
$M'$. We show that 
\beq
\label{fixtrace7}
\fix \p = \fix \p'.
\eeq
Let $u = a_1 \ldots a_k \in
\fix\p$ with $a_1, \ldots, a_k \in A$. Then 
$$a_1\ldots a_k = u = u\p^p = (a_1\p^p)\ldots (a_k\p^p)$$
and $1 \notin A\p$ yields $\{ a_1, \ldots, a_k\} = \{ a_1\p^p, \ldots,
a_k\p^p \}$. By (\ref{fixtrace4}), we get $a_1, \ldots, a_k \in A_0$
and so $u \in \fix\p'$. The converse inclusion holds trivially, hence
$\fix \p = \fix \p'$. 

If $A_0 = A$, then $\fix \p$ is finitely generated by Case I, hence we
may assume that $A_0 \subset A$. By the induction hypothesis, $\fix
\p'$ is finitely generated, and so is $\fix \p$.

\medskip
\noindent
\underline{Case III}: $1 \in A\p$.

\medskip
\noindent
Write 
$$A_1 = A \cap 1\p\inv, \quad A_2 = A \setminus A_1.$$
Let $M''$ be the (trace) submonoid of $M$ generated by $A_2$. Write
$\pi = \pi_{A_2}$ and let $\p'' = (\p\pi)|_{M''}$. Clearly, $\p''
\in \endo M''$ and $A_1 
\neq \emptyset$ implies $|A_2| < |A|$. By the induction hypothesis, $\fix
\p''$ is finitely generated. 
We claim that
\beq
\label{fixtrace6}
\fix\p = (\fix\p'')\p.
\eeq
Let $u \in \fix\p$. We may factor $u = u_0a_1u_1 \ldots a_ku_k$ with
$a_1, \ldots, a_k \in A_2$ and $u_0, \ldots, u_k \in A_1^*$. It
follows that $u = u\p = (a_1a_2 \ldots a_k)\p$. Now $a_1a_2 \ldots a_k
\in M''$ and 
$$(a_1a_2 \ldots a_k)\p'' = (a_1a_2 \ldots a_k)\p\pi =
u\pi = (u_0a_1u_1 \ldots a_ku_k)\pi = a_1a_2 \ldots a_k,$$
hence $a_1a_2 \ldots a_k \in \fix\p''$ and so $u = (a_1a_2 \ldots
a_k)\p \in (\fix\p'')\p$. Thus $\fix\p \subseteq (\fix\p'')\p$.

Conversely, let $v = a_1a_2 \ldots a_k \in \fix\p''$, with $a_1,
\ldots, a_k \in A_2$. By checking directly on the generators, we get
\beq
\label{fixtrace8}
\pi\p = \p.
\eeq
Hence $v = v\p'' = v\p\pi$ yields $(v\p)\p = v\p\pi\p = v\p$ and
so $v\p \in \fix \p$. Thus $(\fix\p'')\p \subseteq \fix\p$ and so
(\ref{fixtrace6}) holds. Since  $\fix
\p''$ is finitely generated, then also  $\fix
\p$ is finitely generated in this third and last case.
\qed

This proof can be adapted to the case of periodic points:

\bt
\label{pertrace}
Let $(A,I)$ be a finite independence alphabet. Then 
{\rm Per}$\,\p$ is finitely generated for every $\p \in$
  {\rm End}$\, M(A,I)$.
\et

\proof
Write $M = M(A,I)$ and let $m = |A|!$. We show that
\beq
\label{pertrace1}
\per\p = \fix\p^m
\eeq
by induction on $|A|$. Then the claim follows from Theorem \ref{fixtrace}.

The case $|A| = 0$ being trivial, assume
that $|A| > 0$ and the result holds for smaller alphabets. We consider
two cases.

\medskip
\noindent
\underline{Case I}: $1 \notin A\p$.

\medskip
\noindent
We keep the notation introduced in Case II of the proof of Theorem
\ref{fixtrace}. We may assume that $A_0 \subset A$, otherwise 
$\p|_A$ would be a permutation, and since the order of $\p|_A$ must divide
the order of the symmetric group on $A$, which is $m$, we would get
$(\p|_A)^m = 1$ and therefore $\p^m = 1$, yielding $\fix\p^m = M = \per\p$.

If we
replace $\p$ by $\p^n$, then $A_0$ remains the same in view of $\per\p
= \per\p^n$ and so does $M'$. On the other hand, it follows from
(\ref{fixtrace5}) that $\p^n|_{M'} =
(\p|_{M'})^n = (\p')^n$, hence 
\beq
\label{pertrace2}
\fix\p^n = \fix(\p')^n
\eeq
by
applying (\ref{fixtrace7}) to $\p^n$. By the induction
hypothesis, we have $\per\p' = \fix(\p')^{|A_0|!}$. Since
$|A_0|!$ divides $m$, we get 
$$\per\p' = \fix(\p')^{|A_0|!} \subseteq \fix(\p')^m \subseteq
\per\p'$$
and so $\per\p' = \fix(\p')^m$. Together with (\ref{pertrace2}), this yields
$$\per\p = \cup_{n\geq 1}\, \fix\p^n = \cup_{n\geq 1}\, \fix(\p')^n =
\per\p' = \fix(\p')^m = \fix\p^m$$
as required.

\medskip
\noindent
\underline{Case II}: $1 \in A\p$.

\medskip
\noindent
We keep the notation introduced in Case III of the proof of Theorem
\ref{fixtrace}.

Let $u \in \per\p$, say $u \in \fix\p^n$. We may factor $u = u_0a_1u_1
\ldots a_ku_k$ with 
$a_1, \ldots, a_k \in A_2$ and $u_0, \ldots, u_k \in A_1^*$. It
follows that $u = u\p^n = (a_1a_2 \ldots a_k)\p^n$. Now $a_1a_2 \ldots a_k
\in M''$ and (\ref{fixtrace8}) yields $a_1a_2 \ldots a_k = (a_1a_2
\ldots a_k)\p^n\pi = (a_1a_2 \ldots a_k)(\p\pi)^n$ and consequently
$a_1a_2 \ldots a_k \in 
\fix(\p'')^n \subseteq \per\p''$. By the induction hypothesis, we get
$a_1a_2 \ldots a_k \in \fix(\p'')^m$ and so $a_1a_2 \ldots a_k = (a_1a_2
\ldots a_k)\p^m\pi$ in view of (\ref{fixtrace8}). Hence $u\p^m =  v_0a_1v_1
\ldots a_kv_k$ for some $v_0, \ldots, v_k \in A_1^*$. Thus $$u\p^{2m} =
u\p^m\pi\p^m = (v_0a_1v_1
\ldots a_kv_k)\pi\p^m = (u_0a_1u_1
\ldots a_ku_k)\pi\p^m = u\p^m$$ which together with $u\p^n = u$ yields
$$u = u\p^n = u\p^{2n} = \ldots = u\p^{nm} = u\p^{(n-1)m} = \ldots =
u\p^m.$$
Therefore $\per\p = \fix\p^m$ and so (\ref{pertrace1}) holds in all
cases as required.
\qed

\bc
\label{constr}
Let $(A,I)$ be a finite independence alphabet and let $\p \in$
  {\rm End}$\, M(A,I)$. Then we can effectively
compute finite sets of generators for {\rm Fix}$\,\p$  and {\rm Per}$\,\p$.
\ec

\proof
For $\fix\p$, it suffices to remark that all the
morphisms, subsets and submonoids appearing in the induction proof of
Theorem \ref{fixtrace} can be
effectively computed, namely in connection with the key equalities 
(\ref{fixtrace1}), (\ref{fixtrace7}) and (\ref{fixtrace6}). The
periodic case follows from the fixed point case and (\ref{pertrace1}).
\qed

\section{Extension to real traces}

In the late eighties, two ultrametrics were introduced for trace
monoids. One of them, defined by Bonizzoni, Mauri and Pighizzini
\cite{BMP}, is known as the {\em FNF metric}. Given
$u,v \in M(A,I)$, let $u = w_{B_1}\ldots w_{B_m}$ and $v =
w_{C_1}\ldots w_{C_n}$ denote their FNFs. We define
$$r(u,v) = \left\{
\begin{array}{ll}
\max\{ k \geq 0 \mid B_1 = C_1, \ldots, B_k = C_k \}&\mbox{ if }u \neq v\\
\infty&\mbox{ otherwise}
\end{array}
\right.$$
The ultrametric $d$ is defined by 
$d(u,v) = 2^{-r(u,v)}$, using the convention $2^{-\infty} = 0$.

The other metric, defined by Kwiatkowska \cite{Kwi}, is known as the {\em
  prefix metric}. Given $u,v \in M(A,I)$, we say that $u$ is a prefix
of $v$ and write $u \leq_p v$ if $v = uw$ for some $w \in M(A,I)$. 
For every $n \in \N$, denote by $\pref_n(v)$ the set of all prefixes of
$u$ of length $n$. We
define, for all $u,v \in M(A,I)$,  
$$r'(u,v) = \sup\{ n \in \N \mid \pref_n(u) = \pref_n(v) \}$$
and $d'(u,v) = 2^{-r'(u,v)}$.
It is well known that, for a finite dependence alphabet, these metrics
are uniformly equivalent (i.e. the 
identity mappings between $(M(A,I),d)$ and $(M(A,I),d')$ are uniformly
continuous), see e.g. \cite{KK}. 
We recall that a mapping $\p:(X_1,d_1) \to (X_2,d_2)$ between metric spaces
is {\em uniformly continuous} if
$$\forall \varepsilon > 0\; \exists \delta > 0 \; \forall x,y \in X_1\;
(d_1(x,y) < \delta \Rw d_2(x\p,y\p) < \varepsilon).$$
An important particular case is given by {\em contractions},
i.e. $d_2(x\p,y\p) \leq d_1(x,y)$ for all $x,y \in X_1$.

In \cite{Kwi}, Kwiatkowska also showed that the completion of
$(M(A,I),d')$ is obtained by adding the {\em infinite real traces} and
is a compact metric space. By standard topology results, this
completion is up to homeomorphism the completion of
$(M(A,I),d)$, hence we shall describe it with respect to the FNF
metric, favoured in this paper.

Let $\partial M(A,I)$ consist of all infinite sequences of the form
$w_{B_1} w_{B_2} \ldots$ such
  that:
\bi
\item
$B_i$ is a clique for every $i \geq 1$;
\item
for all $i \geq 2$ and $a \in B_i$, there exists some $b \in
B_{i-1}$ such that $(a,b) \notin I$.
\ei
Note that $w_{B_1} \ldots w_{B_n}$ is a Foata
normal form for every $n \geq 1$. We say that $\hat{M}(A,I) = M(A,I)
\cup \partial M(A,I)$ is the set of {\em real traces} on $(A,I)$. 

The metric $d$ extends to $\hat{M}(A,I)$ in the obvious way, and it is
easy to check that $(\hat{M}(A,I),d)$ is complete: given a Cauchy
sequence $(u^{[n]})_n$ with $u^{[n]} = w_{B_{n1}} w_{B_{n2}} \ldots$,
it follows easily that each sequence $(w_{B_{nk}})_k$ is stationary
with limit, say, $w_{B_k}$, and we get
$w_{B_1} w_{B_2} \ldots = \lim_{n\to\infty} u^{[n]}$. Since $w_{B_1}
w_{B_2} \ldots = \lim_{n\to\infty} w_{B_1} \ldots w_{B_n}$, then
$(\hat{M}(A,I),d)$ is indeed the completion of $(M(A,I),d)$ and
therefore compact by the aforementioned Kwiatkowska's result. We may
refer to $\partial M(A,I)$ as the {\em boundary} of $M(A,I)$.

It is well known that $M(A,I)$
acts continuously on the left of $\partial M(A,I)$. Indeed, it follows
easily from standard topological arguments which we sketch here. We
start by noting that 
\beq
\label{coac}
d(uv,uw) \leq d(v,w)
\eeq
holds afor all $u,v,w \in M(A,I)$ (i.e. multiplication by a constant
on the left constitutes a contraction of $(M(A,I),d)$). Since
$(\hat{M}(A,I),d)$ is the completion of $(M(A,I),d)$, it follows easily
that the action 
$$\begin{array}{rcl}
M(A,I) \times \partial M(A,I)&\to&\partial M(A,I)\\
(u,\lim_{n\to\infty} v_n)&\mapsto&\lim_{n\to\infty} \, uv_n
\end{array}$$
is well defined, and in view of (\ref{coac}) it turns out to be
continuous. We shall call this action {\em mixed product}.

We remark that the mixed product is left cancellative: if $uX = uX'$
with $u \in M(A,I)$ and $X,X' \in \hat{M}(A,I)$,
then $X = X'$. Indeed, we can extend (\ref{quinta}) to real traces the
obvious way, and left cancellativity becomes a simple exercise.

Kwiatowska's approach, on the other hand, leads to a geometric
description of the boundary reminiscent of 
the theory of hyperbolic groups:

Given a monoid $M$ generated by a subset $A$, the {\em Cayley graph}
$\cay_AM$ has vertex set $M$ and labelled directed edges of the form
$m \mapright{a} ma$ for all $m \in M$ and $a \in A$. Let $\pi:A^* \to
M$ be the canonical epimorphism. We say that a
word $u \in A^*$ is a {\em geodesic} of $\cay_AM$ if $1 \mapright{u}
u\pi$ has shortest length among the
paths connecting $1$ to $u\pi$ in $\cay_AM$. Note that geodesics $1
\mapright{} u\pi$ need not to be unique! An infinite word $\alpha \in
A^{\omega}$ is a {\em ray} of $\cay_AM$ if every finite prefix of $\alpha$
is a geodesic. In the particular case of a trace monoid $M(A,I)$, when
we take the canonical generating set $A$, all relations are
length-preserving and so every finite (respectively infinite) word
represents a geodesic (respectively a ray). The operators $\pref_n$ can
be extended in the obvious way to infinite words and an equivalence
relation can be defined on the set of rays (in our case, $A^{\omega}$)
by
$$\alpha \equiv \beta \hspace{1cm}\mbox{ if $\pref_n(\alpha) =
  \pref_n(\beta)$ for every }n \geq 0.$$
Then the boundary of $M(A,I)$ can be viewed as the set of $\equiv$-equivalence
classes of $A^{\omega}$.

We introduce now a subclass 
of subsets of $\hat{M}(A,I)$ which
generalizes the usual concept of rational subset of a monoid. 
We say that $Y
\subseteq \hat{M}(A,I)$ is
{\em mp-rational} if $Y$ can be obtained from finite subsets of
$\hat{M}(A,I)$ by applying finitely many times the operators union, product,
star and mixed product.
It follows easily that $L$ is 
mp-rational if and only if $L =
L_0 \cup (\cup_{i=1}^n L_iX_i)$ for some rational subsets $L_0,\ldots,L_n$ of
$M(A,I)$ and $X_1,\ldots,X_n \in \partial M(A,I)$. 

For every $u \in M(A,I)$, let $u\xi$ denote the {\em content} of $u$,
i.e. the set of letters occurring in $u$.
We define a symmetric relation $\sim_I$ on $M(A,I)$ by
$$u \sim_I v \hspace{1cm}\mbox{ if }u\xi \times v\xi \subseteq I.$$

\bt
\label{unifcont}
Let $(A,I)$ be a finite independence alphabet and let $\p \in$ {\rm
  End}$\, M(A,I)$. Then the following conditions are equivalent:
\bi
\item[(i)] $\p$ is uniformly continuous with respect to $d$;
\item[(ii)] $\p$ is a contraction with respect to $d$;
\item[(iii)] for all $a,b,c \in A$,
\beq
\label{unifcont1}
(c \leq_p b\p \; \wedge \; c \sim_I a\p) \; \Rw \; (a,b) \in I.
\eeq
\ei
\et

\proof
(i) $\Rw$ (iii). Suppose that (\ref{unifcont1}) fails for some $a,b,c
\in A$. For every $n \in \N$, let $u_n = a^n$ and $v_n = a^nb$. Since
$(a,b) \notin I$, then
$r(u_n,v_n) = n$ and so $d(u_n,v_n) = 2^{-n}$. On the other hand, it
follows from $c \leq_p b\p$ and $c \sim_I a\p$ that $c \not\leq_p u_n\p$
and $c \leq_p v_n\p$, hence the FNFs of $u_n\p$ and
$v_n\p$ differ at the first components. It follows that 
$r(u_n\p,v_n\p) = 0$ and so $d(u_n\p,v_n\p) = 1$. Therefore $\p$ is
not uniformly continuous.

(iii) $\Rw$ (ii). Assume that (\ref{unifcont1}) holds. We say that an
occurrence of a letter in a word $u \in A^*$ has {\em height} $k$ if
it ends up in the $k$th component when we compute the FNF of $u$ by
applying the relations from (\ref{sexta}). Note that we are not
allowed to swap consecutive occurrences of the same letter! 

Let $u =
a_1\ldots a_n$ with $a_1,\ldots,a_n \in A$ and assume that $a_i$ has
height $k$ in $u$. We show that every occurrence of a letter in
$a_i\p$ has height $\geq k$ in $(a_1\p)\ldots (a_n\p)$. We use
induction on $k$. The case $k = 1$ holding trivially, assume that $k >
1$ and the claim holds for $k-1$. There exists some $j < i$ such that
$(a_j,a_i) \notin I$ and $a_j$ has height $k-1$ in $u$. Consider an
occurrence of a letter $c$ in $a_i\p$. Suppose first that this
occurrence has height 1 in $a_i\p$. We may write $c \leq_p
a_i\p$. Since $(a_j,a_i) \notin I$, (\ref{unifcont1}) yields $c
\not\sim_I a_j\p$ and so the height of our occurrence of $c$ is
greater than the height of some occurrence of a letter in $a_j\p$,
which is $\geq k-1$ by the induction hypothesis. On the other hand, if
our occurrence of $c$ has height $> 1$ in $a_i\p$, the first letter in
$a_i\p$ must have height $\geq k$ by the previous case, and we get
height $> k$ for our occurrence. This completes the induction process
and the proof of our claim.

It follows that
\beq
\label{unifcont2}
p \geq k \hspace{.3cm} \Rw \hspace{.3cm} r((w_{B_1}\ldots w_{B_k})\p,
(w_{B_1}\ldots w_{B_p})\p) \geq k 
\eeq
holds whenever $w_{B_1}\ldots w_{B_p}$ is a FNF. Indeed, the freshly proven claim
implies that every occurrence of a letter in $(w_{B_{k+1}}\ldots
w_{B_p})\p$ has height $> k$ in $(w_{B_1}\ldots
w_{B_k})\p(w_{B_{k+1}}\ldots
w_{B_p})\p = (w_{B_1}\ldots w_{B_p})\p)$ and we get (\ref{unifcont2}).

Now it follows easily that
\beq
\label{unifcont3}
r(u\p,v\p) \geq r(u,v)
\eeq
holds for all $u,v \in M(A,I)$. Indeed, if $r(u,v) = k$, we can write $u =
w_{B_1}\ldots w_{B_k}w_{C_1}\ldots w_{C_p}$ and  $v =
w_{B_1}\ldots w_{B_k}w_{D_1}\ldots w_{D_q}$ in FNF. By
(\ref{unifcont2}), we get $r((w_{B_1}\ldots w_{B_k})\p, u\p) \geq k$
and $r((w_{B_1}\ldots w_{B_k})\p, v\p) \geq k$, hence $r(u\p,v\p) \geq
k$ and so (\ref{unifcont2}) holds. Thus $d(u\p,v\p) \leq d(u,v)$ and
so $\p$ is a contraction.

(ii) $\Rw$ (i). Trivial.
\qed

\bc
\label{decuc}
Let $(A,I)$ be a finite independence alphabet and let $\p \in$ {\rm
  End}$\, M(A,I)$. Then the following conditions are decidable:
\bi
\item[(i)] $\p$ is uniformly continuous with respect to $d$;
\item[(ii)] $\p$ is a contraction with respect to $d$.
\ei
\ec

\proof
Immediate from Theorem \ref{unifcont}, since condition (iii) is
obviously decidable.
\qed

Assume now that $\p \in \endo M(A,I)$ is uniformly continuous with
respect to $d$. Since
$(\hat{M}(A,I),d)$ is the completion of $(M(A,I),d)$, $\p$ admits a
unique continuous extension $\Phi$ to $(\hat{M}(A,I),d)$. By
continuity, we must have
\beq
\label{coex}
X\Phi = \lim_{n\to\infty} u_n\p
\eeq
whenever $X \in \partial M(A,I)$ and $(u_n)_n$ is a sequence on
$M(A,I)$ satisfying $X = \lim_{n\to\infty} u_n$.

Given $Y \subseteq M(A,I)$, let $\oo{Y}$ denote the topological
closure of $Y$ in $(\hat{M}(A,I),d)$. It is immediate that $\fix\Phi$
is closed: if $X =\lim_{n \to \infty} X_n$ with every $X_n \in
\fix\Phi$, then $X\Phi = \lim_{n \to \infty} X_n\p = \lim_{n \to \infty}
X_n = X$. It follows that $\oo{\fix\p}
\subseteq \oo{\fix\Phi} = \fix\Phi$.

Note that, for every $u \in M(A,I)$, the sequence $(u^n)_n$ is
Cauchy and therefore convergent. We denote its limit by $u^{\omega}$.

Our study of $\fix\Phi$ starts with the case of free commutative
monoids.

\bl
\label{fcmo}
Let $M$ be a free
  commutative monoid of finite rank and let $\p \in$ {\rm
  End}$\, M$ be uniformly continuous with
respect to $d$. Let $u \in M$ and 
$$Y_{\p,u} = \{ X \in \hat{M} \mid u(X\Phi) = X \}.$$
Then $Y_{\p,u}$ is mp-rational.
\el

\proof
Let $A$ be the basis of $M$. 
We show that
\beq
\label{fcmo2}
Y_{\p,u} = \bigcup_{B \subseteq A} L_Bw_B^{\omega}
\eeq
for some rational subsets $L_B$ of $M\pi_{A\setminus B}$. Since $\hat{M} = \cup_{B
  \subseteq A} (A\setminus B)^*w_B^{\omega}$, it is just the rationality of the
subsets $L_B$ which is at stake. Clearly, a necessary condition for
$L_B \neq \emptyset$ is $B\p\xi = B\xi$, hence we assume the
latter. Assume that $A\setminus B = \{ a_1, \ldots, a_k\}$.
Consider the equation $u((a_1^{x_1}\ldots a_k^{x_k}w_B^{\omega})\Phi) =
a_1^{x_1}\ldots a_k^{x_k}w_B^{\omega}$ on the variables $x_1,\ldots, x_k
\in \mathbb{N}$.
The computation of $L_B$ is done through the system of equations
$$|u|_{a_i} + x_1|a_1\p|_{a_i} + \ldots + x_k|a_k\p|_{a_i} = x_i \quad (i =
1,\ldots,k),$$
and linear diophantine systems such as this are known to have
semilinear (therefore rational) 
sets of solutions \cite{Sim}. Therefore (\ref{fcmo2}) holds and so
$\fix\Phi$ is mp-rational. 
\qed

We say that a graph is of {\em type T} if it has no full
subgraphs of one of the following forms
$$
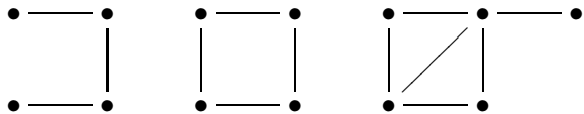
\xymatrix{
\bullet \ar@{-}[r] & \bullet \ar@{-}[d] & \bullet \ar@{-}[r] & \bullet 
& \bullet \ar@{-}[r] \ar@{-}[d] & \bullet \ar@{-}[r] \ar@{-}[d] &
\bullet \\
\bullet \ar@{-}[r] & \bullet & \bullet \ar@{-}[r] \ar@{-}[u] & \bullet
\ar@{-}[u] 
& \bullet \ar@{-}[ur] \ar@{-}[r] & \bullet &
}$$

Write $\Delta_A = \{ (a,a) \mid a \in A \}$. The next result gives a
complete solution for the case when $\Gamma(A,I)$ is of type T.

\bt
\label{infi}
Let $(A,I)$ be a finite independence alphabet such that $\Gamma(A,I)$
is of type T. Then the following 
conditions are equivalent: 
\bi
\item[(i)] for every $\p \in$ {\rm
  End}$\, M(A,I)$, uniformly continuous with
respect to $d$,
there exists some mp-rational $Y \subseteq
\hat{M}(A,I)$ such that
{\rm Fix}$\,\Phi = \oo{{\rm Fix}\,\p} \cup Y$;
\item[(ii)] $I \cup \Delta_A$ is transitive;
\item[(iii)] $\Gamma(A,I)$ is a disjoint union of complete graphs;  
\item[(iv)] $M(A,I)$ is a free product of finitely many free
  commutative monoids of finite rank. 
\ei
\et

\proof
(i) $\Rw$ (ii). Write $M = M(A,I)$.
Suppose that $I \cup \Delta_A$ is not transitive. Then there exist
distinct $a,b,c \in A$ such that $(a,b),(b,c) \in I$ but $(a,c) \notin
I$. Let $\p \in \endo M$ be defined by 
$$x\p = \left\{
\begin{array}{ll}
xb&\mbox{ if }(x,b) \in I\mbox{ and }x \not\sim_I (ac)\\
x&\mbox{ otherwise }
\end{array}
\right.$$
for $x \in A$.

To show that $\p$ is well defined, i.e. $(x,y) \in I \Rw (xy)\p =
(yx)\p$, the only nontrivial case occurs when $x\p = bx$ and $y\p =
y$. Suppose that $by \neq yb$ in $M$. We have $(x,b) \in I$ and from
$x \not\sim_I (ac)$ we may
assume without loss of generality that $xa \neq ax$ in $M$. Hence $a
\edge b \edge x \edge y$ is a 4-vertex subgraph of $\Gamma(A,I)$ and
from $by \neq yb$ and $xa \neq ax$ it follows that the full subgraph
induced by $\{ a,b,x,y \}$ has one of the first two forbidden
configurations, a contradiction. Thus $by = yb$ in $M$ and it follows
easily that $\p$ is well defined.

To show that $\p$ is uniformly continuous, let $x,y,z \in A$ satisfy
$z \leq_p x\p$ and $z \sim_I y\p$. We must show that $(x,y) \in
I$. Suppose not. The case $z = x$ leads to immediate contradiction,
hence we may assume that $x\p = xb$ 
and $z = b$. Obviously, we may assume that $y\p = y$, and also $xa
\neq ax$ in $M$. Now $b \sim_I y\p = y$ yields $y \sim_I (ac)$ and so
$c \neq x$. It follows that $\Gamma(A,I)$ has a subgraph of the form
$$
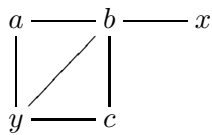
\xymatrix{
a \ar@{-}[r] \ar@{-}[d] & b \ar@{-}[r] \ar@{-}[d] & x \\
y \ar@{-}[ur] \ar@{-}[r] & c & 
}$$
where all the vertices must be
distinct. This cannot be a full subgraph since it is of the third
forbidden type. But the only potential edge that has not been excluded yet is $x
\edge c$, when $a \edge y \edge c \edge x$ would be a full
subgraph of the first forbidden type. Hence we reach a
contradiction in both cases and so  $(x,y) \in
I$. Thus $\p$ satisfies
(\ref{unifcont1}) and is therefore uniformly continuous.

It is immediate that $\fix\p$ is generated by a subset of  $A\setminus
\{ a,c\}$, hence none of the fixed points 
$V_n = ((ab)^n(bc)^n)^{\omega}$
belongs to $\oo{\fix\p}$. 

Suppose that $\fix\Phi = \oo{\fix\p} \; \cup \; L_0 \; \cup \; L_1X_1
\; \cup \; \ldots \; 
\cup \; L_mX_m$ for some rational subsets $L_0, \ldots,L_m$ of $M$ and
$X_1, \ldots, X_m \in \partial M$. Then there exist distinct fixed points
$V_n$ and $V_k$ belonging to the same subset $L_iX_i$. But then both
$V_n$ and $V_k$ must share (infinite) suffixes with $X_i$ and
therefore with each other by transitivity, a contradiction since $n \neq
k$. Therefore condition (i) must fail.

(ii) $\Rw$ (i). Assume that $I \cup \Delta_A$ is transitive and let $\p
\in \endo M$ be uniformly continuous. In view of Lemma \ref{fcmo}, we
may assume that $M$ is noncommutative. We may also assume that $\p$ is
nontrivial. With these assumptions, we claim that $1 \notin
A\p$. Indeed, suppose that $a\p = 1$ for some $a \in A$. Let $b \in
A$. By
(\ref{unifcont1}), we must have $b\p = 1$ whenever $(a,b) \notin
I$. Since $M$ is noncommutative, $\Gamma(A,I)$ has at least two
connected components, and so $c\p = 1$ for every vertex $c$ which is not in the
connected component of $a$. But now, replacing $a$ by one of these
$c$, we must have $d\p = 1$ for every vertex $d$ in the
connected component of $a$ and so $\p$ would be trivial. Thus $1 \notin
A\p$.

Let $A = A_1 \cup \ldots \cup A_r$ be the decomposition of the vertex
set of $\Gamma(A,I)$ in its connected components, and let $M_j$ denote
the free commutative monoid on $A_j$. Let 
$$J = \{ j \in \{ 1,\ldots,r\} \mid M_j\p \subseteq M_j \}.$$
For every $j \in J$, let $\p_j = \p|_{M_j}$. Then $\p_j$ is a
(uniformly continuous) endomorphism of $M_j$, and its continuous
extension $\Phi_j$ to the completion $\hat{M}_j \subseteq \hat{M}$ is
a restriction of $\Phi$. 

Let $P = J \times \{ 0,\ldots, 2^r-1\}$.
For every $(j,k) \in P$, define
$$\begin{array}{ll}
C_{jk} = \displaystyle\bigcup_{i=1}^r \{& B \subseteq A_i \mid w_B\p =
w_Bz_B \mbox{ for some }z_B 
\neq 1\\ 
&\mbox{and }z_B\p^k \in M_j \}
\end{array}$$
and let
$$D = \displaystyle\bigcup_{i = 1}^r \{ B \subseteq A_i \mid \exists
\lim_{n\to\infty} w_B\p^n \}.$$ 

We prove that
\beq
\label{infi4}
\begin{array}{lll}
\fix\Phi&=&\oo{\fix\p} \; \cup (\displaystyle\bigcup_{(j,k) \in P}
\bigcup_{B \in C_{jk}} 
  (\fix\p)w_Bz_B(z_B\p)\ldots (z_B\p^{k-1}) Y_{\p_j,z_B\p^k})\\
&\cup&((\fix\p)\{ \lim_{n\to\infty} 
w_B\p^n \mid B \in D \}).
\end{array}
\eeq
The opposite inclusion follows from straightforward checking, we
consider only the case $X = uw_Bz_B(z_B\p)\ldots (z_B\p^{k-1})X'$ with
$(j,k) \in P$, $B \in C_{jk}$, $u \in \fix\p$ and $X' \in Y_{\p_j,z_B\p^k}$. 
Then
$$X\Phi = (u\p)w_Bz_B(z_B\p)\ldots (z_B\p^{k-1})(z_B\p^{k})(X'\Phi) =
uw_Bz_B(z_B\p)\ldots (z_B\p^{k-1})X' = X$$
as required.

Now take
$X \in \fix\Phi$ and write $X = w_{B_1}w_{B_2}\ldots $ in FNF. Without
loss of generality, we may assume that $w_{B_1} \notin \fix\p$. 
By continuity and
(\ref{unifcont2}), we get $r(w_{B_1}\p,
(w_{B_1}w_{B_2}\ldots)\Phi) \geq 1$ and so $r(w_{B_1}\p,
w_{B_1}w_{B_2}\ldots) \geq 1$. Hence we may write $w_{B_1}\p =
w_{B_1}v$ with $v \neq 1$.

For every $u \in M$, let $u\oo{\xi} = \{ i \in \{ 1,\ldots,r \} \mid
u\xi \cap A_i \neq \emptyset \}$. We claim that
\beq
\label{infi5}
u\oo{\xi} = v\oo{\xi} \; \Rw \; u\p\oo{\xi} = v\p\oo{\xi}
\eeq
holds for all $u,v \in M$. Indeed, suppose that $a \in A_i$ occurs in
$u$. Then some $b \in A_i$ must occur in
$v$. Since $(a\p)(b\p) = (b\p)(a\p)$, it follows easily (directly or
by using Levi's Lemma \cite{DR}) that $a\p\oo{\xi} = b\p\oo{\xi}$. Hence
$u\p\oo{\xi} \subseteq v\p\oo{\xi}$ and (\ref{infi5}) follows from
symmetry. 

Next we define a directed graph $\Omega$ having as vertices the nonempty
subsets of $\{ 1,\ldots,r \}$ and edges $R \mapright{} S$ whenever
$u\oo{\xi} = R$ implies $u\p\oo{\xi} = S$. By (\ref{infi5}), and since
$1 \notin A\p$, $\Omega$ is well defined. Note that each vertex of $\Omega$ has
outdegree 1. 

Suppose first that there exists some $k \geq 0$ and $j \in \{
1,\ldots,r\}$ such that
\beq
\label{flb2}
\{ j \} = v\p^k\oo{\xi} = v\p^{k+1}\oo{\xi} = \ldots
\eeq
It follows that $j
\in J$. Moreover, we may assume
that $k < 2^r$ because after following $2^r-1$ edges in $\Omega$ we
must be in a cycle, and starting at $v\oo{\xi}$ the cycle is
necessarily trivial. Hence $(j,k) \in P$, $B \in C_{jk}$ and we may
write $z_{B_1} = v$. Now, for every $\ell \geq k$ we get
$$w_{B_1}\p^{\ell} = (w_{B_1}z_{B_1})\p^{\ell-1} =
(w_{B_1}z_{B_1}(z_{B_1}\p))\p^{\ell-2} = \ldots =
w_{B_1}z_{B_1}(z_{B_1}\p) \ldots (z_{B_1}\p^{\ell-1}),$$
and so
\beq
\label{infi6}
w_{B_1}w_{B_2}\ldots = (w_{B_1}w_{B_2}\ldots)\Phi^{\ell} = 
w_{B_1}z_{B_1}(z_{B_1}\p) \ldots
(z_{B_1}\p^{\ell-1})((w_{B_2}w_{B_3}\ldots)\Phi^{\ell}).
\eeq
Hence
$$\begin{array}{l}
w_{B_1}z_{B_1}(z_{B_1}\p) \ldots
(z_{B_1}\p^{\ell-1})((w_{B_2}w_{B_3}\ldots)\Phi^{\ell}) =
w_{B_1}w_{B_2}\ldots\\
\hspace{1.5cm}= w_{B_1}z_{B_1}(z_{B_1}\p) \ldots
(z_{B_1}\p^{\ell})((w_{B_2}w_{B_3}\ldots)\Phi^{\ell+1}),
\end{array}$$
yielding by left cancellativity 
\beq
\label{flb}
(w_{B_2}w_{B_3}\ldots)\Phi^{\ell} =
(z_{B_1}\p^{\ell})((w_{B_2}w_{B_3}\ldots)\Phi^{\ell+1}).
\eeq
Iterating (\ref{flb}) for $\ell = k, k+1, \ldots, k'$, we get
$$(w_{B_2}w_{B_3}\ldots)\Phi^k = (z_{B_1}\p^k)(z_{B_1}\p^{k+1})\ldots
(z_{B_1}\p^{k'})((w_{B_2}w_{B_3}\ldots)\Phi^{k'+1}).$$
Since $1 \notin A\p$, it follows easily that 
$$(w_{B_2}w_{B_3}\ldots)\Phi^k =
(z_{B_1}\p^k)(z_{B_1}\p^{k+1})\ldots$$
and so (\ref{flb2}) yields $(w_{B_2}w_{B_3}\ldots)\Phi^k \in
M_j$. Considering (\ref{flb}) for $\ell = k$, we get 
$(w_{B_2}w_{B_3}\ldots)\Phi^k \in Y_{\p_j,z_{B_1}\p^k}$. Together with
(\ref{infi6}), this implies that $X$ belongs to the right
hand side of (\ref{infi4}).

Thus we may assume that the sequence $(v\p^n\oo{\xi})_n$ never
stabilizes on a singular set. For every $n \geq 1$, write $v_n = w_{B_1}\p^n =
w_{B_1}v(v\p) \ldots (v\p^{n-1})$. It follows that the number of
alternating connected components in the sequence
$(v_n)_n$ increases unboundedly. Since $v_n \leq v_{n+1}$ for every
$n$, it follows easily that $(v_n)_n$ is a Cauchy sequence and
therefore convergent in $\hat{M}$. Thus $B_1 \in D$.
Now
$$X = w_{B_1}w_{B_2}\ldots = (w_{B_1}w_{B_2}\ldots)\Phi^n =
v_n((w_{B_2}w_{B_3}\ldots)\Phi^n)$$
yields $X = \lim_{n\to \infty} v_n((w_{B_2}w_{B_3}\ldots)\Phi^n)$. 
Since the number of
alternating connected components in $(v_n)_n$ increases unboundedly,
we immediately get
$X = \lim_{n\to \infty} v_n = \lim_{n\to \infty} w_{B_1}\p^n$ and so
(\ref{infi4}) holds.

In view of Theorem \ref{fixtrace} and Lemma \ref{fcmo}, it follows
that $\fix\Phi$ is the union of $\oo{\fix\p}$ with an mp-rational
subset of $\hat{M}$.

(ii) $\iff$ (iii) $\iff$ (iv). Immediate.
\qed

\bc
\label{fifi}
Let $M$ be a free product of finitely many
  commutative monoids of finite rank and let $\p \in$ {\rm
  End}$\, M$ be uniformly continuous with
respect to $d$. If {\rm Fix}$\,\p$ is finite, so is {\rm Fix}$\,\Phi$.
\ec

\proof
First we note that if $\fix\p$ is finite then $\fix\p = \{ 1 \}$ and
so $\oo{\fix\p} = \{ 1 \}$.
In view of (\ref{infi4}), it suffices to show that every $Y_{\p_j,u}$ in
the conditions of Lemma \ref{fcmo} is finite. 
Write $Y_{\p_j,u} = \bigcup_{B \subseteq A_j} L_Bw_B^{\omega}$
as in (\ref{fcmo2}) and suppose that $L_B$ is infinite for some $B
\subseteq A_j$. Assume that $A_j\setminus B = \{ a_1, \ldots, a_k\}$.
Then the proof of Lemma \ref{fcmo} shows that the system of equations
$$|u|_{a_i} + x_1|a_1\p_j|_{a_i} + \ldots + x_k|a_k\p_j|_{a_i} = x_i \quad (i =
1,\ldots,k)$$
has infinitely many solutions $(x_1,\ldots,x_k ) \in
\mathbb{N}^k$. By Dickson's Lemma \cite{Dic}, there exist some
distinct solutions $(x_1,\ldots,x_k), (y_1,\ldots,y_k)$ such that $x_1
\geq y_1, \ldots, x_k \geq y_k$. It follows that 
$$(x_1-y_1)|a_1\p'|_{a_i} + \ldots + (x_k-y_k)|a_k\p'|_{a_i} = x_i-y_i \quad (i =
1,\ldots,k)$$
and so by the proof of Lemma \ref{fcmo} we get $\fix\p_j \neq \{ 1
\}$ and so $\fix\p \neq \{ 1
\}$, a contradiction. Thus
$Y_{\p_j,u}$ is always finite and so is $\fix\Phi$.
\qed

\section*{Acknowledgements}

The authors acknowledge support from the European Regional Development
Fund through the 
programme COMPETE and by the Portuguese Government through the FCT --
Funda\c c\~ao para a Ci\^encia e a Tecnologia under the project
PEst-C/MAT/UI0144/2011. 
The first author also acknowledges the support
of the FCT project SFRH/BPD/65428/2009.

\end{document}